\documentclass[10pt,reqno]{amsart}

\usepackage{amsmath,amsthm}
\usepackage{amssymb,txfonts} % or use "`pxfonts"' instead of "`txfonts"'
\usepackage{euscript}

\numberwithin{equation}{subsection}

\newcommand{\sgsp}{\renewcommand{\baselinestretch}{1}\tiny\normalsize}
\newcommand{\sqsp}{\renewcommand{\baselinestretch}{1.5}\tiny\normalsize}

\newcommand{\relphantom[1]}{\mathrel{\phantom{#1}}}
\newtheorem{thm}[subsection]{Theorem}

\newtheorem{cor}[subsection]{Corollary}

\newtheorem*{deligne}{Deligne's Conjecture}

\theoremstyle{definition}

\newcommand{\cat}[1]{{\EuScript #1}}
\newcommand{\cA}{\cat{A}}

\newcommand{\cC}{\cat{C}}
\newcommand{\cH}{\cat{H}}
\newcommand{\cO}{\cat{O}}

\newcommand{\cF}{\cat{F}}

\newcommand{\cS}{\cat{S}}

\newcommand{\bF}{\mathbf{F}}

\newcommand{\CH}{\mathrm{CH}_{\mathrm{Hopf}}}
\newcommand{\HH}{HH_{\mathrm{Hopf}}}
\newcommand{\dCH}{d_{\mathrm{CH}}}

\DeclareMathOperator{\Id}{Id}
\DeclareMathOperator{\Hom}{Hom}

\DeclareMathOperator{\CB}{CB}

\begin{document}
\title{Brace operations and Deligne's Conjecture for module-algebras}
\author{Donald Yau}
%\keywords{Module-algebra, Gerstenhaber algebra, brace algebra, homotopy $G$-algebra, Deligne's conjecture}

\begin{abstract}
It is observed that Kaygun's Hopf-Hochschild cochain complex for a module-algebra is a brace algebra with multiplication.  As a result, (i) an analogue of Deligne's Conjecture holds for module-algebras, and (ii) the Hopf-Hochschild cohomology of a module-algebra has a Gerstenhaber algebra structure.
\end{abstract}

%\subjclass[2000]{}
\email{dyau@math.ohio-state.edu}
\address{Department of Mathematics, The Ohio State University Newark, 1179 University Drive, Newark, OH 43055, USA}

%\date{\today}
\maketitle
\sqsp

%%================%%
%%                %%
%%  Introduction  %%
%%                %%
%%================%%

\section{Introduction}
\label{sec:intro}

Let $H$ be a bialgebra and let $A$ be an associative algebra.  The algebra $A$ is said to be an $H$-module-algebra if there is an $H$-module structure on $A$ such that the multiplication on $A$ becomes an $H$-module morphism.  For example, if $S$ denotes the Landweber-Novikov algebra \cite{landweber,novikov}, then the complex cobordism $MU^*(X)$ of a topological space $X$ is an $S$-module-algebra.  Likewise, the singular mod $p$ cohomology $H^*(X; \bF_p)$ of a topological space $X$ is an $\cA_p$-module-algebra, where $\cA_p$ denotes the Steenrod algebra associated to the prime $p$ \cite{es,milnor}.  Other similar examples from algebraic topology can be found in \cite{boardman}.  Important examples of module-algebras from Lie and Hopf algebras theory can be found in, e.g., \cite[V.6]{kassel}.

In \cite{kaygun}, Kaygun defined a Hochschild-like cochain complex $\CH^*(A,A)$ associated to an $H$-module-algebra $A$, called the \emph{Hopf-Hochschild cochain complex}, that takes into account the $H$-linearity.  In particular, if $H$ is the ground field, then Kaygun's Hopf-Hochschild cochain complex reduces to the usual Hochschild cochain complex $C^*(A,A)$ of $A$ \cite{hochschild}.  Kaygun \cite{kaygun} showed that the Hopf-Hochschild cohomology of $A$ shares many properties with the usual Hochschild cohomology.  For example, it can be described in terms of derived functors, and it satisfies Morita invariance.

The usual Hochschild cochain complex $C^*(A,A)$ has a very rich structure.  Namely, it is a brace algebra with multiplication \cite{gv}.  Combined with a result of McClure and Smith \cite{ms} concerning the singular chain operad associated to the little squares operad $\cC_2$, the brace algebra with multiplication structure on $C^*(A,A)$ leads to a positive solution of Deligne's Conjecture \cite{deligne}.  Also, passing to cohomology, the brace algebra with multiplication structure implies that the Hochschild cohomology modules $HH^*(A,A)$ form a Gerstenhaber algebra, which is a graded version of a Poisson algebra.  This fact was first observed by Gerstenhaber \cite{ger}.

The purpose of this note is to observe that Kaygun's Hopf-Hochschild cochain complex $\CH^*(A,A)$ of a module-algebra $A$ also admits the structure of a brace algebra with multiplication.  As in the classical case, this leads to a version of Deligne's Conjecture for module-algebras.  Also, the Hopf-Hochschild cohomology modules $\HH^*(A,A)$ form a Gerstenhaber algebra.  When the bialgebra $H$ is the ground field, these structures reduce to the ones in Hochschild cohomology.

A couple of remarks are in order.  First, there is another cochain complex $\cF^*(A)$ that can be associated to an $H$-module-algebra $A$ \cite{yau}.  The cochain complex $\cF^*(A)$ is a differential graded algebra.  Moreover, it controls the deformations of $A$, in the sense of Gerstenhaber \cite{ger2}, with respect to the $H$-module structure, leaving the algebra structure on $A$ fixed.  It is not yet known whether $\cF^*(A)$ is a brace algebra with multiplication and whether the cohomology modules of $\cF^*(A)$ form a Gerstenhaber algebra.

Second, the results and arguments here can be adapted to module-coalgebras, comodule-algebras, and comodule-coalgebras.  To do that, one replaces the crossed product algebra $X$ (\S \ref{subsec:cpa}) associated to an $H$-module-algebra $A$ by a suitable crossed product (co)algebra \cite{ak1,ak2,ak3} and replaces Kaygun's Hopf-Hochschild cochain complex by a suitable variant.

\subsection{Organization}
\label{subsec:organization}

The rest of this paper is organized as follows.

In the following section, we recall the construction of the Hopf-Hochschild cochain complex $\CH^*(A,A)$ from Kaygun \cite{kaygun}.  In Section \ref{sec:operad}, it is observed that $\CH^*(A,A)$ has the structure of an operad with multiplication (Theorem \ref{thm:operad with mult}).  This leads in Section \ref{sec:brace} to the desired brace algebra with multiplication structure on $\CH^*(A,A)$ (Corollary \ref{cor:brace mult}).  Explicit formulas for the brace operations are given.

In Section \ref{sec:HG}, it is observed that the brace algebra with multiplication structure on $\CH^*(A,A)$ leads to a homotopy $G$-algebra structure (Corollary \ref{cor:HG}).  The differential from this homotopy $G$-algebra and the Hopf-Hochschild differential are then identified, up to a sign (Theorem \ref{thm:diff}).  Passing to cohomology, this leads in Section \ref{sec:G} to a Gerstenhaber algebra structure on the Hopf-Hochschild cohomology modules $\HH^*(A,A)$ (Corollary \ref{cor:G}).  The graded associative product and the graded Lie bracket on $\HH^*(A,A)$ are explicitly described.

In the final section, by combining our results with a result of McClure and Smith \cite{ms}, a version of Deligne's Conjecture for module-algebras is obtained (Corollary \ref{cor:deligne}).  This section can be read immediately after Section \ref{sec:brace} and is independent of Sections \ref{sec:HG} and \ref{sec:G}.

%%===================%%
%%                   %%
%%  Hopf-Hochschild  %%
%%                   %%
%%===================%%

\section{Hopf-Hochschild cohomology}
\label{sec:Hopf}

In this section, we fix some notations and recall from \cite[Section 3]{kaygun} the Hopf-Hochschild cochain complex associated to a module-algebra.

%%%%%%%%%%%%%%%%%%%%%%%%
\subsection{Notations}
\label{subsec:notations}

Fix a ground field $K$ once and for all.  Tensor product and vector space are all meant over $K$.

Let $H = (H, \mu_H, \Delta_H)$ denote a $K$-bialgebra with associative multiplication $\mu_H$ and coassociative comultiplication $\Delta_H$.  It is assumed to be unital and counital, with its unit and counit denoted by $1_H$ and $\varepsilon \colon H \to K$, respectively.

Let $A = (A, \mu_A)$ denote an associative, unital $K$-algebra with unit $1_A$ (or simply $1$).

In a coalgebra $(C, \Delta)$, we use Sweedler's notation \cite{sweedler} for comultiplication: $\Delta(x) = \sum x_{(1)} \otimes x_{(2)}$, $\Delta^2(x) = \sum x_{(1)} \otimes x_{(2)} \otimes x_{(3)}$, etc.

These notations will be used throughout the rest of this paper.

%%%%%%%%%%%%%%%%%%%%%%%%%%%
\subsection{Module-algebra}
\label{subsec:mod-alg}

Recall that the algebra $A$ is said to be an \emph{$H$-module-algebra} \cite{d,kassel,mont,sweedler} if and only if there exists an $H$-module structure on $A$ such that $\mu_A$ is an $H$-module morphism, i.e.,
   \begin{equation}
   \label{eq:mod-alg axiom}
   b(a_1 a_2) \,=\, \sum (b_{(1)}a_1)(b_{(2)}a_2)
   \end{equation}
for $b \in H$ and $a_1, a_2 \in A$.  It is also assumed that $b(1_A) = \varepsilon(b) 1_A$ for $b \in H$.

We will assume that $A$ is an $H$-module-algebra for the rest of this paper.

%%%%%%%%%%%%%%%%%%%%%%%%%%%%%%%%%%%%
\subsection{Crossed product algebra}
\label{subsec:cpa}

Let $X$ be the vector space $A \otimes A \otimes H$.  Define a multiplication on $X$ \cite[Definition 3.1]{kaygun} by setting
   \begin{equation}
   \label{eq:crossed product alg}
   (a_1 \otimes a_1^\prime \otimes b^1)(a_2 \otimes a_2^\prime \otimes b^2) \,\buildrel \text{def} \over=\,
   \sum a_1 \left(b^1_{(1)}a_2\right) \otimes \left(b^1_{(3)} a_2^\prime\right)a_1^\prime \otimes b^1_{(2)}b^2
   \end{equation}
for $a_1 \otimes a_1^\prime \otimes b^1$ and $a_2 \otimes a_2^\prime \otimes b^2$ in $X$.  It is shown in \cite[Lemma 3.2]{kaygun} that $X$ is an associative, unital $K$-algebra, called the \emph{crossed product algebra}.

Note that if $H = K$ ($=$ the trivial group bialgebra $K \lbrack \lbrace e \rbrace \rbrack$), then $X$ is just the enveloping algebra $A \otimes A^{\mathrm{op}}$, where $A^{\mathrm{op}}$ is the opposite algebra of $A$.

The algebra $A$ is an $X$-module via the action
   \[
   (a \otimes a^\prime \otimes b) a_0 \,=\, a(ba_0)a^\prime
   \]
for $a \otimes a^\prime \otimes b \in X$ and $a_0 \in A$.  Likewise, the vector space $A^{\otimes (n+2)}$ is an $X$-module via the action
   \[
   (a \otimes a^\prime \otimes b)(a_0 \otimes \cdots \otimes a_{n+1}) \,=\,
   \sum ab_{(1)}a_0 \otimes b_{(2)}a_1 \otimes \cdots \otimes b_{(n+1)}a_n \otimes b_{(n+2)} a_{n+1}a^\prime
   \]
for $a_0 \otimes \cdots \otimes a_{n+1} \in A^{\otimes (n+2)}$.

%%%%%%%%%%%%%%%%%%%%%%%%
\subsection{Bar complex}
\label{subsec:bar}

Consider the chain complex $\CB_*(A)$ of vector spaces with $\CB_n(A) = A^{\otimes (n+2)}$, whose differential $d_n^{\mathrm{CB}} \colon \CB_n(A) \to \CB_{n-1}(A)$ is defined as the alternating sum $d_n^{\mathrm{CB}} = \sum_{j=0}^n (-1)^j \partial_j$, where
   \[
   \partial_j(a_0 \otimes \cdots \otimes a_{n+1}) \,=\,
   a_0 \otimes \cdots \otimes (a_ja_{j+1}) \otimes \cdots \otimes a_{n+1}.
   \]

It is mentioned above that each vector space $\CB_n(A) = A^{\otimes (n+2)}$ is an $X$-module.  Using the module-algebra condition \eqref{eq:mod-alg axiom}, it is straightforward to see that each $\partial_j$ is $X$-linear.  Therefore, $\CB_*(A)$ can be regarded as a chain complex of $X$-modules.

Note that in the case $H = K$, the chain complex $\CB_*(A)$ of $A \otimes A^{\mathrm{op}}$-modules is the usual bar complex of $A$.

%%%%%%%%%%%%%%%%%%%%%%%%%%%%%%%%%%%%%%%%%%%%
\subsection{Hopf-Hochschild cochain complex}
\label{subsec:HH}

The \emph{Hopf-Hochschild cochain complex of $A$ with coefficients in $A$} is the cochain complex of vector spaces:
   \begin{equation}
   \label{eq:CH}
   (\CH^*(A,A), \dCH) \,\buildrel \text{def} \over=\, \Hom_X((\CB_*(A), d^{\mathrm{CB}}), A).
   \end{equation}
Its $n$th cohomology module, denoted by $\HH^n(A,A)$, is called the \emph{$n$th Hopf-Hochschild cohomology of $A$ with coefficients in $A$}.

When $H = K$, the cochain complex $(\CH^*(A,A), \dCH)$ is the usual Hochschild cochain complex of $A$ with coefficients in itself \cite{hochschild}, and $\HH^n(A,A)$ is the usual Hochschild cohomology module.

In what follows, we will use the notation $\CH^*(A,A)$ to denote (i) the Hopf-Hochschild cochain complex $(\CH^*(A,A), \dCH)$, (ii) the sequence $\lbrace \CH^n(A,A) \rbrace$ of vector spaces, or (iii) the graded vector space $\oplus_n \CH^n(A,A)$.  It should be clear from the context what $\CH^*(A,A)$ means.

%%==========%%
%%          %%
%%  Operad  %%
%%          %%
%%==========%%

\section{Algebraic operad}
\label{sec:operad}

The purpose of this section is to show that the vector spaces $\CH^*(A,A)$ in the Hopf-Hochschild cochain complex of an $H$-module-algebra $A$ with self coefficients has the structure of an operad with multiplication.

%%%%%%%%%%%%%%%%%%%%%%%
\subsection{Operads}
\label{subsec:operads}

Recall from \cite{may1,may2} that an \emph{operad} $\cO = \lbrace \cO(n), \gamma, \Id \rbrace$ consists of a sequence of vector spaces $\cO(n)$ $(n \geq 1)$ together with structure maps
   \begin{equation}
   \label{eq:gamma}
   \gamma \colon \cO(k) \otimes \cO(n_1) \otimes \cdots \otimes \cO(n_k) \,\to\, \cO(n_1 + \cdots + n_k),
   \end{equation}
for $k, n_1, \ldots , n_k \geq 1$ and an \emph{identity element} $\Id \in \cO(1)$, satisfying the following two axioms.
   \begin{enumerate}
   \item The structure maps $\gamma$ are required to be \emph{associative}, in the sense that
   \begin{equation}
   \label{eq:operad1}
   \begin{split}
   \gamma(\gamma(f; g_{1,k}); h_{1,N})
   & =\, \gamma(f; \gamma(g_1; h_{1,N_1}), \ldots , \\
   &\relphantom{} \relphantom{}  \relphantom{}  \relphantom{}  \gamma(g_i; h_{N_{i-1}+1,N_i}), \ldots , \gamma(g_k; h_{N_{k-1}+1,N_k})).
   \end{split}
   \end{equation}
Here $f \in \cO(k)$, $g_i \in \cO(n_i)$, $N = n_1 + \cdots + n_k$, and $N_i = n_1 + \cdots + n_i$.  Given elements $x_i, x_{i+1}, \ldots$, the symbol $x_{i,j}$ is the abbreviation for the sequence $x_i, x_{i+1}, \ldots , x_j$ or $x_i \otimes \cdots \otimes x_j$ whenever $i \leq j$.
   \item The identity element $\Id \in \cO(1)$ is required to satisfy the condition that the linear map
   \begin{equation}
   \label{eq:operad2}
   \gamma(-; \Id, \ldots , \Id) \colon \cO(k) \to \cO(k)
   \end{equation}
is equal to the identity map on $\cO(k)$ for each $k \geq 1$.
   \end{enumerate}

What is defined above is usually called a \emph{non-$\Sigma$ operad} in the literature.

%%%%%%%%%%%%%%%%%%%%%%%%%%%%%%%%%%%%%%%
\subsection{Operad with multiplication}
\label{subsec:operad with m}

Let $\cO$ be an operad.  A \emph{multiplication} on $\cO$ \cite[Section 1.2]{gv} is an element $m \in \cO(2)$ that satisfies
   \begin{equation}
   \label{eq:mult}
   \gamma(m; m, \Id) \,=\, \gamma(m; \Id, m).
   \end{equation}
In this case, $(\cO, m)$ is called an \emph{operad with multiplication}.

%%%%%%%%%%%%%%%%%%%%%%%%%%%%%%%%%%%%%%%%%%%%%%%%%%%%%%%%%%%%%%%%%
\subsection{Operad with multiplication structure on $\CH^*(A,A)$}
\label{subsec:operad structure}

In what follows, in order to simplify the typography, we will sometimes write $\cC(n)$ for the vector space $\CH^n(A,A)$.  To show that the vector spaces $\CH^*(A,A)$ form an operad with multiplication, we first define the structure maps, the identity element, and the multiplication.

\begin{description}
\item[Structure maps] 
For $k, n_1, \ldots , n_k \geq 1$, define a map
   \begin{equation}
   \label{eq:structure map}
   \gamma \colon \cC(k) \otimes \cC(n_1) \otimes \cdots \otimes \cC(n_k) \to \cC(N)
   \end{equation}
by setting
   \begin{multline}
   \label{eq:structure def}
   \gamma(f; g_{1,k})(a_{0,N+1}) \\
   \buildrel \text{def} \over=\,
   f\left(a_0 \otimes g_1(1 \otimes a_{1,n_1} \otimes 1) \otimes \cdots \otimes g_i(1 \otimes a_{N_{i-1}+1, N_i} \otimes 1) \otimes \cdots \otimes a_{N+1}\right).
   \end{multline}
Here the notations are as in the definition of an operad above, and each $a_i \in A$.

\item[Identity element] 
Let $\Id \in \cC(1)$ be the element such that
   \begin{equation}
   \label{eq:identity element}
   \Id(a_0 \otimes a_1 \otimes a_2) \,=\, a_0 a_1 a_2.
   \end{equation}
This is indeed an element of $\cC(1)$, since the identity map on $A$ is $H$-linear.

\item[Multiplication]
Let $\pi \in \cC(2)$ be the element such that
   \begin{equation}
   \label{eq:mult element}
   \pi(a_0 \otimes a_1 \otimes a_2 \otimes a_3) \,=\, a_0 a_1 a_2 a_3.
   \end{equation}
This is indeed an element of $\cC(2)$, since the multiplication map $A^{\otimes 2} \to A$ on $A$ is $H$-linear.
\end{description}

\begin{thm}
\label{thm:operad with mult}
The data $\cC = \lbrace \cC(n), \gamma, \Id \rbrace$ forms an operad.  Moreover, $\pi \in \cC(2)$ is a multiplication on the operad $\cC$.
\end{thm}

\begin{proof}
It is immediate from \eqref{eq:structure def} and \eqref{eq:identity element} that $\gamma(-; \Id^{\otimes k})$ is the identity map on $\cC(k)$ for each $k \geq 1$.

To prove associativity of $\gamma$, we use the notations in the definition of an operad and compute as follows:
   \begin{equation}
   \begin{split}
   \gamma&(\gamma(f; g_{1,k}); h_{1,N})(a_0 \otimes \cdots \otimes a_{M+1}) \\
   & =\, \gamma(f; g_{1,k})(a_0 \otimes \cdots \otimes h_j(1 \otimes a_{M_{j-1}+1,M_j} \otimes 1) \otimes \cdots \otimes a_{M+1}) \\
   & =\, f(a_0 \otimes \cdots \otimes g_i(1 \otimes z_i \otimes 1) \otimes \cdots \otimes a_{M+1}) \\
   & =\, \gamma(f; \ldots , \gamma(g_i; h_{N_{i-1}+1,N_i}), \ldots)(a_0 \otimes \cdots \otimes a_{M+1}).
   \end{split}
   \end{equation}
Here the element $z_i$ $(1 \leq i \leq k)$ is given by
   \begin{equation}
   \label{eq:zi}
   \begin{split}
   z_i
   & \,=\, \bigotimes_{l=N_{i-1}+1}^{N_i} h_l(1 \otimes a_{M_{l-1}+1,M_l} \otimes 1) \\
   & \,=\, h_{N_{i-1}+1}(1 \otimes a_{M_{N_{i-1}}+1,M_{N_{i-1}+1}} \otimes 1) \otimes \cdots \otimes h_{N_i}(1 \otimes a_{M_{N_i-1}+1,M_{N_i}} \otimes 1).
   \end{split}
   \end{equation}
This shows that $\gamma$ is associative and that $\cC = \lbrace \cC(n), \gamma, \Id \rbrace$ is an operad.

To see that $\pi \in \cC(2)$ is a multiplication on $\cC$, one observes that both $\gamma(\pi; \pi, \Id)(a_0 \otimes \cdots \otimes a_4)$ and $\gamma(\pi; \Id, \pi)(a_0 \otimes \cdots \otimes a_4)$ are equal to the product $a_0 a_1 a_2 a_3 a_4$.

This finishes the proof of Theorem \ref{thm:operad with mult}.
\end{proof}

%%=========%%
%%         %%
%%  Brace  %%
%%         %%
%%=========%%

\section{Brace algebra}
\label{sec:brace}

The purpose of this section is to show that the graded vector space $\CH^*(A,A)$ admits the structure of a brace algebra with multiplication.

%%%%%%%%%%%%%%%%%%%%%%%%%%%
\subsection{Brace algebra}
\label{subsec:brace}

For a graded vector space $V = \oplus_{n=1}^\infty V^n$ and an element $x \in V^n$, set $\deg x = n$ and $\vert x \vert = n - 1$.  Elements in $V^n$ are said to have \emph{degree $n$}.

Recall from \cite[Definition 1]{gv} that a \emph{brace algebra} is a graded vector space $V = \oplus V^n$ together with a collection of brace operations $x \lbrace x_1, \ldots , x_n \rbrace$ of degree $-n$, satisfying the associativity axiom:
   \[
   x\lbrace x_{1,m} \rbrace \lbrace y_{1,n} \rbrace
   =\, \sum_{0 \leq i_1 \leq \cdots \leq i_m \leq n} (-1)^\varepsilon x \lbrace y_{1,{i_1}}, x_1 \lbrace y_{{i_1 + 1},{j_1}} \rbrace, y_{j_1 + 1}, \\
   \ldots , y_{i_m}, x_m\lbrace y_{{i_m + 1},{j_m}} \rbrace, y_{{j_m + 1},n} \rbrace.
   \]
Here the sign is given by $\varepsilon = \sum_{p=1}^m \left(\vert x_p \vert \sum_{q=1}^{i_p} \vert y_q \vert\right)$.

%%%%%%%%%%%%%%%%%%%%%%%%%%%%%%%%%%%%%%%%%%%%%%
\subsection{Brace algebra with multiplication}
\label{subsec:brace mult}

Let $V = \oplus V^n$ be a brace algebra.  A \emph{multiplication} on $V$ \cite[Section 1.2]{gv} is an element $m \in V^2$ such that
   \begin{equation}
   \label{eq:mm=0}
   m \lbrace m \rbrace \,=\, 0.
   \end{equation}
In this case, we call $V = (V, m)$ a \emph{brace algebra with multiplication}.

%%%%%%%%%%%%%%%%%%%%%%%%%%%%%%%%%%%%%%%
\subsection{Brace algebra from operad}
\label{susbec:brace from operad}

Suppose that $\cO = \lbrace \cO(n), \gamma, \Id\rbrace$ is an operad.  Define the following operations on the graded vector space $\cO = \oplus \cO(n)$:
   \begin{equation}
   \label{eq:operadbrace}
   x \lbrace x_1, \ldots , x_n \rbrace
   \,\buildrel \text{def} \over=\,
   \sum (-1)^\varepsilon \gamma(x; \Id, \ldots , \Id, x_1, \Id, \ldots , \Id, x_n, \Id, \ldots , \Id).
   \end{equation}
Here the sum runs over all possible substitutions of $x_1, \ldots , x_n$ into $\gamma(x; \ldots )$ in the given order. The sign is determined by $\varepsilon = \sum_{p=1}^n \vert x_p \vert i_p$,
where $i_p$ is the total number of inputs in front of $x_p$.  Note that
   \[
   \deg x \lbrace x_1, \ldots , x_n \rbrace \,=\, \deg x - n + \sum_{p=1}^n \deg x_p,
   \]
so the operation \eqref{eq:operadbrace} is of degree $-n$.

Proposition 1 in \cite{gv} establishes that the operations \eqref{eq:operadbrace} make the graded vector space $\oplus \cO(n)$ into a brace algebra.  Moreover, a multiplication on the operad $\cO$ in the sense of \S \ref{subsec:operad with m} is equivalent to a multiplication on the brace algebra $\oplus \cO(n)$.  In fact, for an element $m \in \cO(2)$, one has that
     \begin{equation}
     \label{eq:mult brace}
     m \lbrace m \rbrace \,=\, \gamma(m; m, \Id) - \gamma(m; \Id, m).
     \end{equation}
It follows that the condition \eqref{eq:mult} is equivalent to \eqref{eq:mm=0}.  In other words, an operad with multiplication $(\cO, m)$ gives rise to a brace algebra with multiplication $(\oplus \cO(n), m)$.  Combining this discussion with Theorem \ref{thm:operad with mult}, we obtain the following result.

\begin{cor}
\label{cor:brace mult}
The graded vector space $\CH^*(A,A)$ is a brace algebra with brace operations as in \eqref{eq:operadbrace} and multiplication $\pi$ \eqref{eq:mult element}.
\end{cor}

The brace operations on $\CH^*(A,A)$ can be described more explicitly as follows.  For $f \in \cC(k)$ and $g_i \in \cC(m_i)$ $(1 \leq i \leq n)$, we have
   \begin{equation}
   \label{eq:brace explicit}
   f \lbrace g_1, \ldots , g_n \rbrace
   \,=\, \sum (-1)^\varepsilon \gamma(f; \Id^{r_1}, g_1, \Id^{r_2}, g_2, \ldots , \Id^{r_n}, g_n, \Id^{r_{n+1}}),
   \end{equation}
where $\Id^r = \Id^{\otimes r}$.  Here the $r_j$ are given by
   \begin{equation}
   \label{eq:rj}
   r_j \,=\,
   \begin{cases}
   i_1 & \text{ if } j = 1, \\
   i_j - i_{j-1} - 1 & \text{ if } 2 \leq j \leq n , \\
   k - i_n - 1 & \text{ if } j = n+1,
   \end{cases}
   \end{equation}
and
   \begin{equation}
   \label{eq:epsilon}
   \varepsilon \,=\, \sum_{p=1}^n (m_p - 1)i_p.
   \end{equation}
Write $M = \sum_{i=1}^n m_i$ and $M_j = \sum_{i=1}^j m_i$.  Then for an element $a_{0,\,k+M-n+1} \in A^{\otimes (k+M-n)}$, we have
   \begin{equation}
   \label{eq:brace exp}
   \begin{split}
   f \lbrace g_{1,n} \rbrace (a_{0,\,k+M-n+1}) 
   & = \sum (-1)^\varepsilon f(a_{0,i_1} \otimes g_1(1 \otimes a_{i_1+1,\, i_1+m_1} \otimes 1) \otimes \cdots \\
   & \relphantom{} \relphantom{} \relphantom{} \relphantom{} \relphantom{} \otimes a_{i_{j-1} + M_{j-1} - (j-1) + 2, \, i_j + M_{j-1} - j + 1} \\
   & \relphantom{} \relphantom{} \relphantom{} \relphantom{} \relphantom{} \otimes g_j(1 \otimes a_{i_j + M_{j-1} - j + 2, \, i_j + M_j - j + 1} \otimes 1) \otimes \cdots \\
   & \relphantom{} \relphantom{} \relphantom{} \relphantom{} \relphantom{} \otimes a_{i_n + M - n + 2, \, k+M-n+1}).
   \end{split}
   \end{equation}

%%================%%
%%                %%
%%  Gerstenhaber  %%
%%                %%
%%================%%

\section{Homotopy Gerstenhaber algebra}
\label{sec:HG}

The purpose of this section is to observe that the brace algebra with multiplication structure on $\CH^*(A,A)$ induces a homotopy Gerstenhaber algebra structure.

%%%%%%%%%%%%%%%%%%%%%%%%%%%%%%%%%%
\subsection{Homotopy $G$-algebra}
\label{subsec:HG}

Recall from \cite[Definition 2]{gv} that a \emph{homotopy $G$-algebra} $(V,d, \cup)$ consists of a brace algebra $V = \oplus V^n$, a degree $+1$ differential $d$, and a degree $0$ associative $\cup$-product that make $V$ into a differential graded algebra, satisfying the following two conditions.
   \begin{enumerate}
   \item The $\cup$-product is required to satisfy the condition
   \[
   (x_1 \cup x_2) \lbrace y_{1,n} \rbrace
   \,=\, \sum_{k=0}^n (-1)^\varepsilon x_1 \lbrace y_{1,k} \rbrace \cup x_2 \lbrace y_{{k+1},n} \rbrace,
   \]
where $\varepsilon = \vert x_2 \vert \sum_{p=1}^k \vert y_p \vert$, for $x_i, y_j \in V$.
   \item The differential is required to satisfy the condition
   \[
   \begin{split}
   & d(x \lbrace x_{1,n+1} \rbrace) - (dx)\lbrace x_{1,n+1}\rbrace \\
   & \relphantom{} - (-1)^{\vert x \vert} \sum_{i=1}^{n+1} (-1)^{\vert x_1 \vert + \cdots + \vert x_{i-1} \vert} x \lbrace x_1, \ldots , dx_i, \ldots , x_{n+1} \rbrace \\
   & =\, (-1)^{\vert x \vert \vert x_1 \vert + 1} x_1 \cup x \lbrace x_{2,n+1} \rbrace \\
   & \relphantom{} + (-1)^{\vert x \vert} \sum_{i=1}^n (-1)^{\vert x_1 \vert + \cdots + \vert x_{i-1} \vert} x \lbrace x_1, \ldots , x_i \cup x_{i+1}, \ldots , x_{n+1} \rbrace \\
   & \relphantom - x \lbrace x_{1,n} \rbrace \cup x_{n+1}.
   \end{split}
   \]
   \end{enumerate}

%%%%%%%%%%%%%%%%%%%%%%%%%%%%%%%%%%%%%%%%%%%%%%%%%%%%%%%%%%%%%%%%%%%%%%%%
\subsection{Homotopy $G$-algebra from brace algebra with multiplication}
\label{subsec:HG brace}

A brace algebra with multiplication $V = (V, m)$ gives rise to a homotopy $G$-algebra $(V, d, \cup)$ \cite[Theorem 3]{gv}, where the $\cup$-product and the differential are defined as:
   \begin{equation}
   \label{eq:hG}
   \begin{split}
   x \cup y & \,\buildrel \text{def}\over=\, (-1)^{\deg x} m \lbrace x, y \rbrace, \\
   d(x)     & \,\buildrel \text{def}\over=\, m \lbrace x \rbrace - (-1)^{\vert x \vert}x \lbrace m \rbrace.
   \end{split}
   \end{equation}
This applies to the brace algebra $\CH^*(A,A)$ with multiplication $\pi$ (Corollary \ref{cor:brace mult}).

\begin{cor}
\label{cor:HG}
For an $H$-module-algebra $A$, $\cC = (\CH^*(A,A), d, \cup)$ is a homotopy $G$-algebra.
\end{cor}

%%%%%%%%%%%%%%%%%%%%%%%%%%%%%%%%%%%%
\subsection{Comparing differentials}
\label{subsec:comparing d}

At this moment, there are two differentials on the graded vector space $\CH^*(A,A)$, namely, the differential $d^n$ \eqref{eq:hG} induced by the multiplication $\pi$ and the Hopf-Hochschild differential $\dCH^n$.  The following result ensures that the cohomology modules defined by these two differentials are the same.

\begin{thm}
\label{thm:diff}
The equality $\dCH^n = (-1)^{n+1}d^n$ holds for each $n$.
\end{thm}

\begin{proof}
Pick $f \in \CH^n(A,A)$.  Then we have
   \begin{equation}
   \label{eq:df}
   \begin{split}
   d^n f & \,=\, \pi \lbrace f \rbrace \,+\, (-1)^n f \lbrace \pi \rbrace \\
   & =\, (-1)^{n-1} \gamma(\pi; \Id, f) \,+\, \gamma(\pi; f, \Id) \\
   & \relphantom{} \relphantom{} +\, (-1)^n \sum_{i=1}^n (-1)^{i-1} \gamma(f; \Id^{i-1}, \pi, \Id^{n-i}). 
   \end{split}
   \end{equation}
It follows that
   \begin{equation}
   \label{eq:sign df}
   (-1)^{n+1}d^nf \,=\, \gamma(\pi; \Id, f) + (-1)^{n+1} \gamma(\pi; f, \Id)
   + \sum_{i=1}^n (-1)^i \gamma(f; \Id^{i-1}, \pi, \Id^{n-i}).
   \end{equation}
Observe that \cite[page 8]{kaygun}
   \begin{equation}
   \label{eq:linearity}
   g(a_{0,n+1}) \,=\, a_0 g(1 \otimes a_{1,n} \otimes 1)a_{n+1}
   \end{equation}
for $g \in \CH^n(A,A)$.   Using \eqref{eq:linearity} and applying the various terms in \eqref{eq:sign df} to an element $a_{0,n+2} \in \CB_{n+1}(A) = A^{\otimes(n+3)}$, we obtain
   \begin{equation}
   \label{eq:dfa}
   \begin{split}
   \gamma(\pi; \Id, f)(a_{0,n+2}) &\,=\, f(a_0a_1 \otimes a_{2,n+2}), \\
   \gamma(\pi; f, \Id)(a_{0,n+2}) &\,=\, f(a_{0,n} \otimes a_{n+1}a_{n+2}), \\
   \gamma(f; \Id^{i-1}, \pi, \Id^{n-i})(a_{0,n+2}) &\,=\, f(a_{0,i-1} \otimes a_ia_{i+1} \otimes a_{i+2,n+2}).
   \end{split}
   \end{equation}
The Theorem now follows immediately from \eqref{eq:sign df} and \eqref{eq:dfa}.
\end{proof}

\begin{cor}
\label{cor:diff}
There is an isomorphism of cochain complexes 
   \[
   (\CH^*(A,A), \dCH) \, \cong \,(\CH^*(A,A), d).  
   \]
Moreover, the cohomology modules on $\CH^*(A,A)$ defined by the differentials $\dCH$ and $d$ are equal.
\end{cor}

%%=============%%
%%             %%
%%  G-algebra  %%
%%             %%
%%=============%%

\section{Gerstenhaber algebra}
\label{sec:G}

The purpose of this section is to observe that the homotopy $G$-algebra structure on $\CH^*(A,A)$ gives rise to a $G$-algebra structure on the Hopf-Hochschild cohomology modules $\HH^*(A,A)$.

%%%%%%%%%%%%%%%%%%%%%%%%%%%%%%%%%
\subsection{Gerstenhaber algebra}
\label{subsec:G}

Recall from \cite[Section 2.2]{gv} that a \emph{$G$-algebra} $(V, \cup, \lbrack -, - \rbrack)$ consists of a graded vector space $V = \oplus V^n$, a degree $0$ associative $\cup$-product, and a degree $-1$ graded Lie bracket
   \[
   \lbrack -, - \rbrack \colon V^m \otimes V^n \to V^{m+n-1},
   \]
satisfying the following two conditions:
   \begin{equation}
   \label{eq:G axioms}
   \begin{split}
   x \cup y & \,=\, (-1)^{\deg x \deg y}y \cup x, \\
   \lbrack x, y \cup z \rbrack &\,=\, \lbrack x, y \rbrack \cup z + (-1)^{\vert x \vert \deg y} y \cup \lbrack x, z \rbrack.
   \end{split}
   \end{equation}
In other words, the $\cup$-product is graded commutative, and the Lie bracket is a graded derivation for the $\cup$-product.  In particular, a $G$-algebra is a graded version of a Poisson algebra.  This algebraic structure was first studied by Gerstenhaber \cite{ger}.

%%%%%%%%%%%%%%%%%%%%%%%%%%%%%%%%%%%%%%%%%%%%%%%%%%
\subsection{$G$-algebra from homotopy $G$-algebra}
\label{subsec:G HG}

If $(V, d, \cup)$ is a homotopy $G$-algebra, one can define a degree $-1$ operation on $V$,
   \begin{equation}
   \label{eq:bracket}
   \lbrack x, y \rbrack \,\buildrel \text{def} \over=\,
   x \lbrace y \rbrace - (-1)^{\vert x \vert \vert y \vert} y \lbrace x \rbrace.
   \end{equation}
Passing to cohomology, $(H^*(V,d), \cup, \lbrack -,- \rbrack)$ becomes a $G$-algebra (\cite{gv} Corollary 5 and its proof).

Combining the previous paragraph with Corollary \ref{cor:HG} and Corollary \ref{cor:diff}, we obtain the following result.

\begin{cor}
\label{cor:G}
The Hopf-Hochschild cohomology modules $\HH^*(A,A)$ of an $H$-module-algebra $A$ admits the structure of a $G$-algebra.
\end{cor}

This $G$-algebra can be described on the cochain level more explicitly as follows.  Pick $\varphi \in \CH^n(A,A)$ and $\psi \in \CH^m(A,A)$.  Then
   \begin{equation}
   \label{eq:cup Lie G}
   \begin{split}
   (\psi \cup \varphi)(a_{0,m+n+1}) 
   & \,=\, (-1)^{m+n-1} \psi(a_{0,m} \otimes 1) \varphi(1 \otimes a_{m+1,m+n+1}), \\
   \lbrack \psi, \varphi \rbrack & \,=\, \psi \lbrace \varphi \rbrace - (-1)^{(m-1)(n-1)} \varphi \lbrace \psi \rbrace,
   \end{split}
   \end{equation}
where, writing $a = a_{0,m+n}$,
   \begin{equation}
   \label{eq:Lie G}
   \begin{split}
   \psi \lbrace \varphi \rbrace (a) 
   & \,=\, \sum_{i=1}^m (-1)^{(i-1)(n-1)} \psi(a_{0,i-1} \otimes \varphi(1 \otimes a_{i,i+n-1} \otimes 1) \otimes a_{i+n,m+n}), \\
   \varphi \lbrace \psi \rbrace (a) 
   & \,=\, \sum_{j=1}^n (-1)^{(j-1)(m-1)} \varphi(a_{0,j-1} \otimes \psi(1 \otimes a_{j,j+m-1} \otimes 1) \otimes a_{j+m,m+n}).
   \end{split}
   \end{equation}
In particular, if $m = n = 1$, then the bracket operation
   \begin{equation}
   \label{eq:bracket 1}
   \lbrack \psi, \varphi \rbrack(a_{0,2}) \,=\, \psi(a_0 \otimes \varphi (1 \otimes a_1 \otimes 1) \otimes a_2) - \varphi(a_0 \otimes \psi(1 \otimes a_1 \otimes 1) \otimes a_2)
   \end{equation}
gives $\HH^1(A,A)$ a Lie algebra structure.  There is another description of this Lie algebra  in terms of (inner) derivations in \cite[Proposition 3.9]{kaygun}.

%%===========%%
%%           %%
%%  Deligne  %%
%%           %%
%%===========%%

\section{Deligne's Conjecture for module-algebras}
\label{sec:deligne}

The purpose of this section is to observe that a version of Deligne's Conjecture holds for the Hopf-Hochschild cochain complex of a module-algebra.  The original Deligne's Conjecture for Hochschild cohomology is as follows.

\begin{deligne}[\cite{deligne}]
The Hochschild cochain complex $C^*(R,R)$ of an associative algebra $R$ is an algebra over a suitable chain model of May's little squares operad $\cC_2$ \cite{may1}.
\end{deligne}

A positive answer to Deligne's conjecture was given by, among others, McClure and Smith \cite[Theorem 1.1]{ms} and Kaufmann \cite[Theorem 4.2.2]{kaufmann}.  There is an operad $\cH$ whose algebras are the brace algebras with multiplication (\S \ref{subsec:brace mult}).  For an associative algebra $R$, the Hochschild cochain complex $C^*(R,R)$ is a brace algebra with multiplication and hence an $\cH$-algebra.  McClure and Smith showed that $\cH$ is quasi-isomorphic to the chain operad $\cS$ obtained from the little squares operad $\cC_2$ by applying the singular chain functor, thereby proving Deligne's Conjecture.

It has been observed that the Hopf-Hochschild cochain complex $\CH^*(A,A)$ is a brace algebra with multiplication (Corollary \ref{cor:brace mult}).  Therefore, we can use the result of McClure and Smith \cite[Theorem 1.1]{ms} to obtain the following version of Deligne's Conjecture for module-algebras.

\begin{cor}[Deligne's Conjecture for module-algebras]
\label{cor:deligne}
The Hopf-Hochschild cochain complex $\CH^*(A,A)$ of an $H$-module-algebra $A$ is an algebra over the McClure-Smith operad $\cH$ that is a chain model for the little squares operad $\cC_2$.
\end{cor}

%%==============%%
%%              %%
%%  References  %%
%%              %%
%%==============%%

\sgsp

\end{document}